\documentclass[12pt]{article}

 \usepackage{mathrsfs,amsfonts}
 \usepackage[dvips]{color}

 \newtheorem{theorem}{Theorem}[section]
 \newtheorem{lemma}[theorem]{Lemma}
 \newtheorem{corol}[theorem]{Corollary}
 \newtheorem{prop}[theorem]{Proposition}
 \newtheorem{example}[theorem]{Example}

 \def\blemma{\begin{lemma}\sl{}\def\elemma{\end{lemma}}}
 \def\bproposition{\begin{prop}\sl{}\def\eproposition{\end{prop}}}
 \def\btheorem{\begin{theorem}\sl{}\def\etheorem{\end{theorem}}}
 
 \def\bexample{\begin{example}\rm{}\def\eexample{\end{example}}}

 \def\beqlb{\begin{eqnarray}}\def\eeqlb{\end{eqnarray}}
 \def\beqnn{\begin{eqnarray*}}\def\eeqnn{\end{eqnarray*}}

 \def\proof{\noindent{\it Proof.~~}}\def\qed{\hfill$\Box$\medskip}

 \def\<{\langle}\def\>{\rangle}

 \def\mcr{\mathscr}\def\mbb{\mathbb}

 \def\ar{\!\!&}

\begin{document}

\centerline{\Large\bf Perturbations of continuous-time Markov chains}

\bigskip\bigskip

\bigskip

\centerline{Pei-Sen Li}

\smallskip

\centerline{School of Mathematical Sciences, Beijing Normal University,}

\centerline{Beijing 100875, China}

\centerline{E-mail: \tt peisenli@mail.bnu.edu.cn}

\bigskip

{\narrower

\noindent\textit{Abstract.} The equivalence of regularity of a $Q$-matrix with its bounded perturbations is proved and a integration by parts formula is established for the associated Feller minimal transition functions.

\smallskip


\par}

\bigskip


\section{Introduction}

\setcounter{equation}{0}

One of the basic questions in studying continuous-time Markov chains is to find the regularity
criterion, i.e., to investigate the conditions under which the given $Q$-matrix is regular, or,
equivalently, the corresponding Feller minimal process is honest in the sense that the corresponding
transition function $P(t)=\{P_{ij}(t); i,j\in\mbb{N}\}$ satisfies $\sum^\infty_{j=0} P_{ij}(t) =1$ for all $i\geq 0$ and $t\geq 0$. Here we assume the chain has state space $\mbb{N}:=\{0,1,2,\ldots\}$. We refer to Anderson (1991) and Chen (2004) for the general theory of continuous-time Markov chains. In this note we show that the regularity property is preserved under a bounded perturbation of the $Q$-matrix. We also establish a integration by parts formula for the corresponding Feller minimal processes without the regularity condition.

Given two $Q$-matrices $R=(r_{ij}; i,j\in\mbb{N})$ and $A= (a_{ij}; i,j\in\mbb{N})$, we call $Q= (q_{ij}; i,j\in\mbb{N}):=R+A$ the \emph{perturbation of $R$ by $A$}. Throughout this note, we assume all $Q$-matrices are stable and conservative.

The main purpose of this note is to prove the following theorems:

\btheorem\label{t1.3} Suppose that $A$ is a bounded $Q$-matrix. Then $Q=R+A$ is regular if and only if  $R$ is regular.
\etheorem

\btheorem\label{t1.4}
Let $Q(t)=\{Q_{ij}(t); i,j\in\mbb{N}\}$ and $R(t)=\{R_{ij}(t); i,j\in\mbb{N}\}$ be the Feller minimal transition functions of $Q$ and $R$, respectively. Then we have the following integration by parts formula
\beqlb\label{e3.2}
\sum_{k\in \mbb{N}}\int^t_0 R_{ik}(s)a_k Q_{kj}(t-s)ds= \sum_{l\in \mathbb{N},m\neq l}\int^t_0R_{il}(t-v)a_{lm}Q_{mj}(v)dv + R_{ij}(t) - Q_{ij}(t).
 \eeqlb
In particular, when $\sum_{k\in \mbb{N}}\int^t_0 R_{ik}(s)a_k Q_{kj}(t-s)ds< \infty$, we can rewrite (\ref{e3.2}) as
 \beqlb\label{e3.4}
Q(t)-R(t)=\int^t_0 R(s)AQ(t-s)ds.
\eeqlb
\etheorem

The perturbation theory of infinitesimal generators has been a very useful tool in the hands of analysts and physicists. A considerable amount of research has been done on the perturbation of linear operators on a Banach space. The effect on a semigroup by adding a linear operator to its infinitesimal generator was studied by Phillips (1952) and Yan (1988). However, these authors did not show the equivalence of the regularity of a $Q$-matrix with its bounded perturbations. The  integration by parts formula~(\ref{e3.2}) was given by Chen (2004, p510) under a stronger condition.
The $Q$-matrix of the branching processes with immigration and/or resurrection introduced in Li and Chen (2006) can be regarded as the perturbations of a given branching $Q$-matrix.

\bexample\label{e1}
Let $R= (r_{ij}; i, j\in \mbb{N})$ be a branching $Q$-matrix given by
\beqnn
r_{ij}=\left\{
\begin{array}{lcl}
i b_{j-i+1} & & {j\geq i-1, i\ge 1};\\
0 & & \mbox{otherwise},
\end{array} \right.
\eeqnn
where
$$ b_j\geq 0\ (j\neq 0), \qquad\sum^\infty_{n\neq 1} b_j= -b_1\geq0.$$
Let
$A= (a_{ij}; i, j\in \mbb{N})$ be a bounded $Q$-matrix given by
\beqnn
a_{ij}=\left\{
\begin{array}{lcl}
c_{j-i+1} & & {j\geq i, i\ge 1};\\
h_{j} & & {j\geq 0, i= 0};\\
0 & & \mbox{otherwise},
\end{array} \right.
\eeqnn
where
\beqnn
\left\{
\begin{array}{lcl}
h_j\geq 0\ (j\neq 0), & & \sum^\infty_{j=1} h_j=-h_0\geq 0;\\
c_j\geq 0\ (j\neq 0), & & \sum^\infty_{j=1} c_j=-c_0\geq 0.
\end{array} \right.
\eeqnn
Then the $Q$-matrix $Q= (q_{ij} ; i, j\in Z+):=R+A$ is called a branching $Q$-matrix
with immigration and resurrection.
The corresponding continuous-time Markov chain is called a branching
process with immigration and resurrection.
Note that the regularity criterion of the branching $Q$-matrix $R$ is given by Harris (1963). Since $A$ is a bounded $Q$-matrix, by Theorem~\ref{t1.3} we see $Q$ is regular if and only $R$ is regular. This simplifies considerably the proof of Theorem 2.1 in Li and Chen (2006).
\eexample
\section{Bounded perturbations}
In this section, we assume $A$ is a bounded $Q$-matrix. We shall prove that the regularity of $R$ and $Q$ are equivalent. Let $\gamma=\sup_i a_i=-\inf_ia_{ii}$.
Let $q'_{ii} =\gamma-a_i> 0$ and $q'_{ij}=q_{ij}$ for $i\neq j$. Let  $a'_{ij}= a_{ij}+\gamma \delta_{ij}> 0$. Then we have $q'_{ik}= a'_{ik}+(1-\delta_{ik})r_{ik}> 0$,

\bproposition\label{pert2.3} The backward Kolmogorov equation of $Q$ is equivalent to the following equation:
 \beqlb\label{per2.2}
Q_{ij}(t)=\sum_{k\in \mbb{N}}\int^t_0 e^{-(r_i+\gamma) (t-s)} q'_{ik} Q_{kj}(s) ds+\delta_{ij}e^{-(r_i+\gamma) t}.
 \eeqlb
\eproposition

\proof Suppose that $Q(t)=\{Q_{ij}(t); i,j\in\mbb{N}\}$ is a solution of the backward Kolmogorov equation $\partial_tQ(t)=QQ(t)$. Then
 $$
\partial_tQ_{ij}(t)+(r_i+\gamma)Q_{ij}(t)=\sum_{k\in \mbb{N}}q'_{ik}Q_{kj}(t).
 $$
Multiplying both sides by the integrating factor $e^{(r_i+\gamma)t}$, we find
 $$
\partial_t (e^{(r_i+\gamma)t}Q_{ij}(t)) = e^{(r_i+\gamma)t}\sum_{k\in \mbb{N}}q'_{ik}Q_{kj}(t).
 $$
Integrating and dividing both sides by $e^{(r_i+\gamma)t}$ give~(\ref{per2.2}).
Conversely, suppose $Q_{ij}(t)$ is a solution of~(\ref{per2.2}). By differentiating both sides of the equation we get the backward Kolmogorov equation $\partial_tQ(t)=QQ(t)$.\qed

Let $Q(t)=\{Q_{ij}(t); i,j\in\mbb{N}\}$ and $R(t)=\{R_{ij}(t); i,j\in\mbb{N}\}$ be the minimal transition functions of $Q$ and $R$, respectively. By the \emph{second successive approximation scheme}; see, e.g., Chen (2004, p64), we see
 \beqlb\label{per2.3}
Q_{ij}(t)=\sum^\infty_{n=0} Q_{ij}^{(n)}(t)\quad \mbox{and}\quad R_{ij}(t)=\sum^\infty_{n=0} R_{ij}^{(n)}(t),
 \eeqlb
 where
 \beqlb\label{per2.4}
R^{(0)}_{ij}(t)=\delta_{ij}e^{-r_i t},\quad R^{(n+1)}_{ij}(t)= \sum_{k\neq i}\int^t_0 e^{-r_i(t-s)} r_{ik}R^{(n)}_{kj}(s)ds
 \eeqlb
and
 \beqlb\label{per2.5}
Q^{(0)}_{ij}(t)=\delta_{ij}e^{-(r_i+\gamma) t},\quad Q^{(n+1)}_{ij}(t)= \sum_{k\in \mbb{N}}\int^t_0 e^{-(r_i+\gamma)(t-s)} q'_{ik}Q^{(n)}_{kj}(s)ds.
 \eeqlb

\blemma\label{pert2.4} For any $n\ge 0$ we have
 \beqlb\label{per2.6}
Q_{ij}^{(n)}(t)= \sum^{n-1}_{p=0}\sum_{l,k\in \mbb{N}}\int^t_0e^{-\gamma(t-s)}R^{(n-p-1)}_{il}(t-s)a'_{lk} Q^{(p)}_{kj}(s)ds+R_{ij}^{(n)}(t)e^{-\gamma t}
 \eeqlb
with $\sum^{-1}_{p=0}=0$ by convention.
\elemma

\proof For $n=0$, we have (\ref{per2.6}) trivially. Suppose that (\ref{per2.6}) holds for $n=0,1,\cdots,m$.
Recall that $q'_{ik}= a'_{ik}+(1-\delta_{ik})r_{ik}$. By the second equality in (\ref{per2.5}) we have
 \beqnn
Q_{ij}^{(m+1)}(t) \ar=\ar \sum_{k\neq i} \int^t_0 e^{-(r_i+\gamma)(t-s)}r_{ik}Q^{(m)}_{kj}(s)ds \cr
 \ar\ar\qquad\quad
+ \sum_{k\in \mbb{N}} \int^t_0 e^{-(r_i+\gamma)(t-s)}a'_{ik}Q^{(m)}_{kj}(s)ds \cr
 \ar=:\ar  I_1+I_2.
 \eeqnn
By (\ref{per2.6}) and (\ref{per2.4}) we have
 \beqnn
I_1\ar=\ar\sum_{k\neq i} \int^t_0 e^{-(r_i+\gamma)(t-s)} r_{ik}\bigg[\sum^{m-1}_{p=0}\sum_{l,r\in\mbb{N}} \int^s_0e^{-\gamma(s-u)}R^{(m-p-1)}_{kl}(s-u)a'_{lr} Q^{(p)}_{rj}(u)du\bigg]ds\cr
 \ar\ar\qquad\qquad
+ \sum_{k\neq i} \int^t_0 e^{-(r_i+\gamma)(t-s)} r_{ik} R_{kj}^{(m)}(s)e^{-\gamma s} ds \cr
 \ar=\ar
\sum^{m-1}_{p=0}\sum_{l,r\in\mbb{N}}\int^t_0e^{-\gamma(t-u)}\bigg[\sum_{k\neq i}\int^t_u e^{-r_i(t-s)}r_{ik}R_{kl}^{(m-p-1)}(s-u)ds\bigg] a'_{lr} Q^{(p)}_{rj}(u)du\cr
\ar\ar\qquad\qquad
+ \sum_{k\neq i} \int^t_0 e^{-(r_i+\gamma)(t-s)} r_{ik}R_{kj}^{(m)}(s)e^{-\gamma s}ds \cr
 \ar=\ar
\sum^{m-1}_{p=0}\sum_{l,r\in\mbb{N}} \int^t_0e^{-\gamma(t-u)}R_{il}^{(m-p)}(t-u) a'_{lr} Q^{(p)}_{rj}(u)du + R_{ij}^{(m+1)}(t)e^{-\gamma t}.
 \eeqnn
On the other hand, using the first equality in (\ref{per2.4}) we obtain
 \beqnn
I_2\ar=\ar\sum_{k\in \mbb{N}} \int^t_0 e^{-\gamma(t-s)}R_{ii}^{(0)}(t-s) a'_{ik}Q^{(m)}_{kj}(s)ds\cr
\ar=\ar
\sum_{l,k\in \mbb{N}} \int^t_0 e^{-\gamma(t-s)}R_{il}^{(0)}(t-s) a'_{lk}Q^{(m)}_{kj}(s)ds.
\eeqnn
Summing up the above expressions of $I_1$ and $I_2$, we see (\ref{per2.6}) also holds when $n=m+1$. That gives the desired result. \qed

 \bproposition\label{pert2.1}
Let $Q(t)=\{Q_{ij}(t); i,j\in\mbb{N}\}$ and $R(t)=\{R_{ij}(t); i,j\in\mbb{N}\}$ be the minimal transition functions of $Q$ and $R$, respectively. Then $Q_{ij}(t)$ is the unique solution of the following equation
\beqlb\label{per2.1}
Q_{ij}(t)= \sum_{l,k\in \mbb{N}}\int^t_0 e^{-\gamma (t-s)}R_{il}(t-s)a'_{lk}Q_{kj}(s)ds+R_{ij}(t)e^{-\gamma t}.
\eeqlb
 \eproposition
 \proof
We first prove the uniqueness of (\ref{per2.1}). Let $\tilde{Q}_{ij}(t)$ be another solution of (\ref{per2.1}). Let $c_{ij}(t)= |Q_{ij}(t)-\tilde{Q}_{ij}(t)|$ and $c_j(t)= \sup_{i} c_{ij}(t).$
Then we have
 $$c_{ij}(t)\leq \sum_{l, k\in \mbb{N}}\int^t_0 R_{il}(t-s)a'_{lk}c_{kj}(s)ds.$$
Taking the supremum we have
 \beqnn
c_j(t)\leq \sup_i \sum_{l,k\in \mbb{N}} \int^t_0 R_{il}(t-s)a'_{lk}c_j(s)ds
 =
 \gamma\int^t_0 c_j(s)ds.
 \eeqnn
Using Gronwall's inequality we have that $c_j(t)=0$. Thus (\ref{per2.1}) has at most one solution.

Next we will show that $Q_{ij}(t)$ satisfies (\ref{per2.1}).
Using (\ref{per2.3}) and $(\ref{per2.6})$ we have
 \beqnn
Q_{ij}(t)=\sum^\infty_{n=0}\sum^{n}_{p=0}\sum_{l,k\in \mbb{N}}\int^t_0e^{-\gamma(t-s)}R^{(n-p)}_{il}(t-s)a'_{lk} Q^{(p)}_{kj}(s)ds+\sum^{\infty}_{n=0}R^{(n)}_{ij}(t)e^{-\gamma t}.
 \eeqnn
Interchanging the order of summation and using (\ref{per2.3}) again we obtain
 \beqnn
Q_{ij}(t)
\ar=\ar
\sum_{l,k\in \mbb{N}}\int^t_0e^{-\gamma(t-s)}\sum^\infty_{n=p}R^{(n-p)}_{il}(t-s)a'_{lk} \sum^\infty_{p=0}Q^{(p)}_{kj}(s)ds+R_{ij}(t)e^{-\gamma t}\cr
 \ar=\ar
\sum_{l,k\in \mbb{N}}\int^t_0 e^{-\gamma (t-s)}R_{il}(t-s)a'_{lk}Q_{kj}(s)ds+R_{ij}(t)e^{-\gamma t}.
 \eeqnn
That completes the proof.
\qed

\noindent\emph{Proof of Theorem~\rm\ref{t1.3}.}~
Summing up both sides of~(\ref{per2.1}) over $j$, we see that $x_i(t):=\sum^\infty_{j=0}Q_{ij}(t)$ is a solution to the following equation:
 \beqlb\label{per2.7}
x_i(t)= \sum_{l,k\in \mbb{N}}\int^t_0 e^{-\gamma (t-s)}R_{il}(t-s)a'_{lk}x_k(s)ds+e^{-\gamma t}\sum^{\infty}_{j=0}R_{ij}(t).
 \eeqlb
Suppose that $R$ is regular. Then we have $\sum^\infty_{j=0}R_{ij}(t)=1$, so $x_i(t)\equiv 1$ is a solution of~(\ref{per2.7}). Let $\tilde{x}_{i}(t)$ be another solution of~(\ref{per2.7}). Set $c_{i}(t)= |x_{i}(t)-\tilde{x}_{i}(t)|$ and $c(t) = \sup_{i} c_{i}(t).$
By (\ref{per2.7}) we obtain
 $$
c_{i}(t)\leq \sum_{l, k\in \mbb{N}}\int^t_0 R_{il}(t-s)a'_{lk}c_{k}(s)ds.
 $$
Taking the supremum we get
 \beqnn
c(t)\ar\leq\ar \sup_i \sum_{l,k\in \mbb{N}} \int^t_0 R_{il}(t-s)a'_{lk}c(s)ds
=
\gamma\int^t_0 c(s)ds.
 \eeqnn
Using Gronwall's inequality we have $c(t)=0$. Then we see $x_i(t)\equiv 1$ is the unique solution to~(\ref{per2.7}).  Hence $Q$ is regular.

Conversely, suppose that $Q$ is regular. Then $x_i(t)= \sum^\infty_{j=0} Q_{ij}(t)=1$. Let $y_i(t) = \sum^\infty_{j=0} R_{ij}(t)$. From~(\ref{per2.7}) we have
 $$
1-e^{-\gamma t}\leq \int^t_0\gamma e^{-\gamma (t-s)}y_i(t-s)ds.
 $$
Then we must have $y_i(t)\equiv 1$, so $R$ is regular.
\qed

\section{Integration by parts formula}
Recall that $R(t)$ and $Q(t)$ are the Feller minimal transition functions of $R$ and $Q$, respectively. By the second successive approximation scheme; see, e.g. Chen (2004, p64) we have
 \beqlb\label{3.10}
Q_{ij}(t)=\sum^\infty_{n=0} Q_{ij}^{(n)}(t),
 \eeqlb
 where
 \beqlb\label{3.17}
Q^{(0)}_{ij}(t)=\delta_{ij}e^{-q_i t},\quad Q^{(n+1)}_{ij}(t)= \sum_{k\neq i}\int^t_0 e^{-q_i(t-s)} q_{ik}Q^{(n)}_{kj}(s)ds.
 \eeqlb
Let $(\Omega,\mcr{F},\mcr{F}_t, \xi_t, P_i)$ be a realization of $(R_{ij}(t))_{t\geq 0}$.
 \blemma\label{l1}
Let $\sigma^t_s$ denote the number of jumps of the trajectory $t\mapsto \xi_t$ on the interval $(s,t]$. Then for $n\geq0$ we have
 \beqlb\label{3.1}
Q_{ij}^{(n)}(t)
 \ar=\ar\nonumber
\sum^{n-1}_{p=0}\sum_{k\in
\mbb{N}, l\neq k}\int^t_0P_i (M^0_{t-s};A_{n-p-1,k}(0, t-s) )a_{kl} Q^{(p)}_{lj}(s)ds\\
\ar\ar\qquad\qquad\qquad\qquad+P_i (M^0_t;A_{n,j}(0,t) )
 \eeqlb
with $\sum^{-1}_{p=0}=0$ by convention, where $A_{n,j}(s,t)= \{\sigma^t_s=n, \xi_t=j\}$ and $M^r_t=e^{-\int^t_r a(\xi_s)ds}$.
 \elemma
 \proof
For $n=0$ we have (\ref{3.1}) trivially. Suppose that (\ref{3.1}) holds for $n=0,1,\cdots,m$.
By (\ref{3.17}) we have
 \beqnn
Q_{ij}^{(m+1)}(t) \ar=\ar \sum_{k\neq i} \int^t_0 e^{-q_i(t-s)}r_{ik}Q^{(m)}_{kj}(s)ds \cr
 \ar\ar
+ \sum_{k\neq i} \int^t_0 e^{-q_i(t-s)}a_{ik}Q^{(m)}_{kj}(s)ds =:  I_1+I_2.
 \eeqnn
  Denote $\tau=\inf\{t\geq 0:\xi_t\neq \xi_{0}\} $. By the
Markov property we have
  \beqnn
I_1\ar=\ar\sum_{k\neq i}\int^t_0 e^{-q_i(t-s)}r_{ik}
\sum^{m-1}_{p=0}\sum_{r\in \mbb{N},l\neq r} \int^s_0
P_k (M^0_{s-v}
;A_{m-1-p,r}(0,s-v) )
a_{rl}Q_{lj}^{(p)}(v)dvds\cr
 \ar\ar\qquad\qquad
+\sum_{k\neq i} \int^t_0 e^{-q_i(t-s)}r_{ik}
P_k (M^0_{s};A_{m,j}(0,s) )ds\cr
 \ar=\ar
\sum^{m-1}_{p=0}\sum_{r\in \mbb{N},l\neq r}\int^t_0\sum_{k\neq i}\int^{t}_v
e^{-q_i(t-s)}r_{ik}P_k (M^0_{s-v}; A_{m-1-p,r}(0,s-v) )
a_{rl}Q_{lj}^{(p)}(v)dsdv\cr
 \ar\ar\qquad\qquad+\sum_{k\neq i}\int^t_0 e^{-a_i(t-s)}r_i^{-1}r_{ik}P_k (M^0_s;A_{m,j}(0,s) )r_ie^{-r_i(t-s)}ds\cr
 \ar=\ar
\sum^{m-1}_{p=0}\sum_{r\in \mbb{N},l\neq r}\int^t_0 \Big[\sum_{k\neq
i}\int^{t-v}_0
e^{-a_is}r_i^{-1}r_{ik}P_k (M^0_{t-v-s}; A_{m-1-p,r}(0,t-v-s) )\cr
 \ar\ar r_ie^{-r_is}ds \Big]a_{rl}Q_{lj}^{(p)}(v)dv+\sum_{k\neq i}\int^t_0
 e^{-a_is}r_i^{-1}r_{ik}P_k(M^0_{t-s}; A_{m,j}(0,t-s))r_ie^{-r_is}ds\cr
 \ar=\ar
\sum^{m-1}_{p=0}\sum_{r\in \mbb{N},l\neq r}\int^t_0 P_i\Big[e^{-a_i\tau}P_{\xi_\tau}(M^0_{t-v-\tau} ; A_{m-1-p,r}(0,t-v-\tau))\Big]a_{rl}Q_{lj}^{(p)}(v)dv\cr
 \ar\ar\qquad\qquad+P_i\bigg[e^{-a_i\tau}P_{\xi_\tau}(M^0_{t-\tau}; A_{m,j}(0,t-\tau))\bigg]\cr
 \ar=\ar
\sum^{m-1}_{p=0}\sum_{r\in \mbb{N},l\neq r}\int^t_0 P_i\Big[e^{-a_i\tau}P_i (M^{\tau}_{t-v} 1_{A_{m-p,r}(\tau,t-v)}|\mcr{F}_\tau )\Big]a_{rl}Q_{lj}^{(p)}(v)dv\cr
 \ar\ar\qquad\qquad+P_i\bigg[e^{-a_i\tau} P_i(M^\tau_t1_{A_{m+1,j}(\tau,t)}|\mcr{F}_\tau)\bigg]\cr
 \ar=\ar
\sum^{m-1}_{p=0}\sum_{r\in \mbb{N},l\neq r}\int^t_0P_i(M^0_{t-v};A_{m-p,r}(0,t-v))a_{rl}Q_{lj
}^{(p)}(v)dv\cr
\ar\ar\qquad\qquad+P_i(M^0_t;A_{m+1,j}(0,t)).
 \eeqnn
On the other hand, we have
\beqlb\label{3.18} I_2\ar=\ar\sum_{l\neq i} \int^t_0 e^{-q_i(t-v)}a_{il}Q^{(m)}_{lj}(v)dv\cr
\ar=\ar
\sum_{l\neq i}\int^t_0P_i\Big(e^{-a_i(t-v)}1_{\{ \sigma^{t-v}_0=0\}}\Big)a_{il}Q_{lj
}^{(m)}(v)dv\cr
\ar=\ar
\sum_{r\in \mbb{N},l\neq r}\int^t_0P_i (M^0_{t-v};A_{0,r}(0,t-v) )a_{rl}Q_{lj
}^{(m)}(v)dv.
\eeqlb
Summing up the above expressions of $I_1$ and $I_2$ we see (\ref{3.1}) also holds when $n=m+1$. That gives the desired result. \qed

\btheorem\label{t3.2} The Feller minimal transition functions $Q(t)$ and $R(t)$ satisfy the following
equation \beqlb\label{e3.3} Q_{ij}(t)\ar=\ar \sum_{k\in \mathbb{N},l\neq
k} \int^t_0P_i( M^0_{t-s}1_{\{\xi_{t-s}=k\}})a_{kl}Q_{lj}(s)ds
+P_i(M^0_t1_{\{\xi_t=j\}}). \eeqlb

\etheorem
\proof Using (\ref{3.10}) and (\ref{3.1}) we have
 \beqnn
Q_{ij}(t)\ar=\ar \sum^\infty_{n=0}\sum^{n-1}_{m=0}\sum_{k\in \mbb{N},l\neq k}\int^t_0P_i(M^0_{t-s};A_{k,n-m-1}(0,t-s))a_{kl} Q^{(m)}_{lj}(s)ds
+P_i(M^0_t1_{\{\xi_t=j\}}).
 \eeqnn
Interchanging the order of summation we see (\ref{e3.3}) holds.
\qed

\emph{Proof of Theorem~{\rm \ref{t1.4}.}}~
By the Markov property of $\{\xi_t: t\geq 0\}$,
 \beqnn
 \ar\ar\sum_{k\in \mbb{N}}\int^t_0 R_{ik}(s)a_k P_k( M^0_{t-s}1_{\{\xi_{t-s}=j\}})ds\cr
 \ar\ar\qquad=
\int^t_0 P_i\Big[a(\xi_s)P_{\xi_s}( M^0_{t-s}1_{\{\xi_{t-s}=j}\})\Big]ds\cr
 \ar\ar\qquad=
\int^t_0 P_i\Big[a(\xi_s)P_i\Big( M^s_t1_{\{\xi_t=j\}}|\mcr{F}_s\Big)\Big]ds\cr
 \ar\ar\qquad=
\int^t_0P_i\Big[a(\xi_s)M^s_t1_{\{\xi_t=j\}}\Big]ds\cr
 \ar\ar\qquad=
P_i\Big[1_{\{\xi_t=j\}}\int^t_0a(\xi_s)e^{-\int^t_sa(\xi_u)du}ds\Big]\cr
 \ar\ar\qquad=
P_i\Big[1_{\{\xi_t=j\}}\Big(1-e^{-\int^t_0a(\xi_u)du}\Big)\Big]\cr
 \ar\ar\qquad=
R_{ij}(t)-P_i\Big(M^0_t1_{\{\xi_t=j\}}\Big).
 \eeqnn
On the other hand, by the Markov property, we have
\beqnn
\ar\ar\sum_{l\in \mathbb{N},m\neq l}\int^t_0R_{il}(t-v)a_{lm}Q_{mj}(v)dv
\cr
\ar\ar\qquad= \sum_{l\in \mathbb{N},m\neq l} P_i\Big[\int^t_0 1_{\{\xi_{t-v}=l\}}a_{lm}Q_{mj}(v)\Big(1-e^{-\int^{t-v}_0a(\xi_u)du}\Big)dv\Big]\cr
\ar\ar\qquad\quad+
\sum_{l\in \mathbb{N},m\neq l} P_i\Big(M^0_{t-v}1_{\{\xi_{t-v=l}\}}\Big)a_{lm}Q_{mj}(v)dv\cr
\ar\ar\qquad= \sum_{l\in \mathbb{N},m\neq l} P_i\Big[\int^t_01_{\{\xi_{t-v}=l\}}a_{lm}Q_{mj}(v)dv\int^{t- v}_0 a(\xi_s)e^{-\int^{t-v}_sa(\xi_u)du}ds\Big]\cr
\ar\ar\qquad\quad+
\sum_{l\in \mathbb{N},m\neq l} P_i\Big(M^0_{t-v}1_{\{\xi_{t-v=l}\}}\Big)a_{lm}Q_{mj}(v)dv\cr
\ar\ar\qquad= \sum_{l\in \mathbb{N},m\neq l} P_i\Big[\int^t_0a(\xi_s)\int^{t-s}_{0}M^s_{t-v}1_{\{\xi_{t-v=l}\}}a_{lm}Q_{mj}(v)dvds\Big]\cr
\ar\ar\qquad\quad+
\sum_{l\in \mathbb{N},m\neq l} P_i\Big(M^0_{t-v}1_{\{\xi_{t-v=l}\}}\Big)a_{lm}Q_{mj}(v)dv\cr
\ar\ar\qquad= \sum_{l\in \mathbb{N},m\neq l}P_i\Big\{\int^t_0 a(\xi_s)P_i\Big[\int^{t-s}_{0}M^s_{t-v}1_{\{\xi_{t-v}=l\}}a_{lm}Q_{mj}(v)dv|\mcr{F}_s\Big]ds\Big\}\cr
\ar\ar\qquad\quad+
\sum_{l\in \mathbb{N},m\neq l} P_i\Big(M^0_{t-v}1_{\{\xi_{t-v=l}\}}\Big)a_{lm}Q_{mj}(v)dv\cr
\ar\ar\qquad= P_i\Big\{\int^t_0 a(\xi_s)\sum_{l\in \mathbb{N},m\neq l} P_{\xi_s}\Big[\int^{t-s}_0M^0_{t-s-v}1_{\{\xi_{t-s-v=l}\}}a_{lm}Q_{mj}(v)dv\Big]ds\Big\}\cr
\ar\ar\qquad\quad+
\sum_{l\in \mathbb{N},m\neq l} P_i\Big(M^0_{t-v}1_{\{\xi_{t-v=l}\}}\Big)a_{lm}Q_{mj}(v)dv\cr
\ar\ar\qquad= \sum_{k\in \mathbb{N}} \int^t_0 R_{ik}(s)a_k\sum_{l\in \mathbb{N},m\neq l}\int^{t-s}_0 P_k(M^0_{t-s-v}1_{\{\xi_{t-s-v}=l\}})a_{lm}Q_{mj}(v)dvds\cr
\ar\ar\qquad\quad+
\sum_{l\in \mathbb{N},m\neq l} P_i\Big(M^0_{t-v}1_{\{\xi_{t-v=l}\}}\Big)a_{lm}Q_{mj}(v)dv.
\eeqnn
By the above two equations and~(\ref{e3.3}) we obtain~(\ref{e3.2}). Suppose that $$\sum_{k\in \mbb{N}}\int^t_0 R_{ik}(s)a_k Q_{kj}(t-s)ds< \infty.$$ Then subtracting it from both sides of~(\ref{e3.2}) yields~(\ref{e3.4}).
\qed

\noindent\textbf{Acknowledgments.}~ The author would like to thank Professors Mu-Fa Chen, Yong-Hua Mao and Yu-Hui Zhang for their helpful comments. Research supported in part by 985 Project, NSFC(No 11131003, 11501531, 11571043), SRFDP(No 20100003110005) and the Fundamental Research Funds for the Central Universities.

\end{document}